\documentclass[12pt]{article}
\usepackage{amsfonts,amssymb,amsmath,theorem}
\newtheorem{Pa}{Paper}[section]
\newtheorem{Tm}[Pa]{{\bf Theorem}}
\newtheorem{La}[Pa]{{\bf Lemma}}
\newtheorem{Cy}[Pa]{{\bf Corollary}}
\newtheorem{Pn}[Pa]{{\bf Proposition}}
{
\theorembodyfont{\normalfont}
\newtheorem{Rk}[Pa]{{\bf Remark}}
\newtheorem{Ee}[Pa]{{\bf Example}}
\newtheorem{Dn}[Pa]{{\bf Definition}}
\newtheorem{I}[Pa]{{\bf}}

}

\newcommand\qed{\ifhmode\unskip\nobreak\fi\quad  
   \ifmmode\square\else\hbox{$\square$}\fi}      

\let\cal=\mathcal

\renewcommand{\H}{\mathcal H}

\newcommand{\B}{\mathcal B}

\newcommand{\V}{\mathcal V}

\newcommand{\al}{\alpha}

\newcommand{\supp}{\operatorname{supp}}
\newcommand{\Prim}{\operatorname{Prim}}
\newcommand{\End}{\operatorname{End}}
\newcommand{\A}{\mathcal A}

\renewcommand{\P}{\mathcal P}
\newcommand{\C}{\mathbb C}
\newcommand{\R}{\mathbb R}

\newcommand{\N}{\mathbb N}

\numberwithin{equation}{section}

\begin{document}
\begin{center}
{\bfseries \Large Crossed product of a $C^*$-algebra by a semigroup of bounded positive linear maps. Interactions.}

\bigskip
{\large \textsc{B. K.  Kwa\'sniewski} }\\
\bigskip
  Institute of Mathematics,  University  in Bialystok,\\
 ul. Akademicka 2, PL-15-267  Bialystok, Poland\\
 e-mail: bartoszk@math.uwb.edu.pl
 \end{center}

\begin{abstract} The paper  presents  a  construction  of  the  crossed  product  of  a
$C^*$-algebra  by a commutative semigroup  of  bounded positive linear maps  generated by partial isometries.
In particular, it generalizes Antonevich,  Bakhtin, Lebedev's crossed product by an endomorphism, and is related
to Exel's interactions. One of the main goals is the Isomorphism Theorem established in the case of actions by endomorphisms.
\end{abstract}

\medbreak

\medbreak

\vspace{5mm}
\tableofcontents

\section{Introduction}\label{intro}
Recently, in \cite{Ant-Bakht-Leb} A.B. Antonevich, V.I. Bakhtin
 and A.V. Lebedev introduced a  new crossed product  of a $C^*$-algebra by an endomorphism
 (for abbreviation we shall call it ABL-crossed product) which in a sense, see \cite{Ant-Bakht-Leb},
  generalizes all the  previous approaches   to constructions of that kind  in the case of a single endomorphism
  \cite{cuntz}, \cite{CK}, \cite{Paschke}, \cite{Exel1}
\cite{Murphy},  \cite{exel},
\cite{kwa}.  Afterwards, see \cite{Kwa-Leb}, the ABL-crossed product  was adapted  to the case
 of actions by  a semigroup $\Gamma^+$ which is a positive cone of a  totally ordered commutative group $\Gamma$.
  It is fundamental that the ABL-construction arose against a background of R. Exel's crossed product \cite{exel},
  which (was dapted to the semigroup context  by N. S. Larsen \cite{Larsen} and requires
 a  unital $C^*$-algebra $\A$, a semigroup homomorphism $\al:\Gamma^+ \to \End (\A)$ where $\End(\A)$ is the set of
 endomorphisms of $\A$ (with composition as a semigroup operation), and  also
  it depends on a choice of \emph{transfer action},  i. e. a semigroup
  homomorphism $L:\Gamma^+\to \mathrm{PosLin}(\A)$ where $\mathrm{PosLin}(\A)$ is the set of all
  linear bounded positive maps on  $\A$, such that
$$
L_x(\al_x(a)b)=aL_x(b), \qquad \textrm{ for all } \,\, a,b\in \A\,\, \textrm{ and }\, x\in \Gamma^+.
$$
In other words Exel's crossed product is a certain $C^*$-algebra  associated to
 the system $(\A,\Gamma^+,\al, L)$ consisting of  four elements (cf. Example {2-adic}),
  whereas the ABL-crossed product \cite{Ant-Bakht-Leb}, \cite{Kwa-Leb}, depends only on
   the triple $(\A,\Gamma^+,\al)$. The price to pay (which eventually is not that high, see
   \cite{Ant-Bakht-Leb}) is that ABL-crossed product is defined only  for a special class of
   \emph{finely representable systems} $(\A,\Gamma^+, \al)$, see  \cite{Kwa-Leb}, \cite{Ant-Bakht-Leb}.
\\
 A link between Exel's and ABL-crossed product is provided by the result of V.I. Bakhtin
 and A.V. Lebedev  \cite{Bakht-Leb}  which being stated in the  semigroup language \cite{Kwa-Leb}
  says that $(\A,\Gamma^+,\al)$ is finely representable if and only if there exists a
  transfer action $L$ for  $(\A,\Gamma^+,\al)$, such that
\begin{equation}
\al_x(L_x(a))=\al_x(1)a\al_x(1),\quad \qquad \textrm{ for all }\,\,  a\in \A, \,\, x\in \Gamma^+ ,
\end{equation}
in which case $L$ is called a \emph{complete transfer action}. It is important that the complete transfer action, if
it exists, is  unique and $\al$ and $L$  determine uniquely  one another via the formulae
$$
L_x(a)=\al_x^{-1}(\al_x(1)a\al_x(1)), \qquad \al_x(a)=L_x^{-1}(L_x(1)a) \qquad a\in \A,\, \, x\in \Gamma^+,
$$
see \cite[Theorem 2.8]{Bakht-Leb}, \cite[Theorem 2.4]{Kwa-Leb}.
\\

Let us note that, although one can not help feeling that in the
above  picture the action $\al$ is somewhat privileged,  there is
no particular reason to single out $\al$ since we have one-to-one
correspondence $\al \longleftrightarrow L$ (in the ABL-context, of
course). This simple observation is a starting point for the
present article. We attempt to clarify here a number of  questions
which  arise  naturally:
\begin{itemize}
\item Why not carry out the ABL-construction starting with $L$
rather than with $\al$? \item Is it necessary for one of the
elements in the pair $(\al,L)$ to act by multiplicative mappings?
\item What happens if we drop this multiplicativity condition,
which of the results concerning ABL-crossed products can be
carried over then?
\end{itemize}
Furthermore, we are not  simply interested in  generalizing
ABL-crossed product.  We also aim at a powerful tool to study
crossed products the so-called  Isomorphism Theorem \cite{Jun},
\cite{Anton_Lebed}, \cite{top-free}, \cite{kwa} which has not been
studied in the ABL-context  yet.
\par
 We have to mention one more important fact. In \cite{exel-inter} a similar dissatisfaction  of an asymmetry between
 actions and transfer actions in  the construction of  Exel's crossed product  led R. Exel
  to an object which he called interaction. Simply, due to the author of  \cite{exel-inter} \emph{interaction}
   is a pair $(\V,\H)$ of two positive bounded linear maps on a $C^*$-algebra $\A$ such that
$$
\V \circ  \H\circ \V= \V,\qquad \H \circ  \V\circ \H= \H
$$
$$
\V|_{\H(\A)}\,\,\, \textrm{ and }\,\,\, \H|_{\V(\A)}\,\,\, \textrm{ are multiplicative}.
$$
It is quite striking that a connection of the article
\cite{exel-inter} with the present paper is completely analogous to
that of Exel's crossed product with ABL-crossed product (which will
become clear during the further reading).

The paper is organized as follows.
\\
In Section \ref{the_interactions} we convert Exel's notion of
interaction to the semigroup case and  present some of its properties.
Then in Section \ref{complette interactions} we define  complete interactions,
 explain its connection with complete transfer actions, and give a few  characteristics of this notion.
Section \ref{endom} is devoted to finely representable actions and
associated crossed products. Here, we define finely
representability of an action $\V$ and then show that it implies
the existence  of (necessarily unique)  action $\H$ such that the
pair $(\V,\H)$ is a complete interaction. We also develop some
terminology and  facts concerning  the internal structure of the
crossed product, which we use later in Section \ref{2.1} to obtain
a necessary and sufficient condition for a representation of the
crossed product to be faithful.  The final Section \ref{topolo
free}  is dedicated to the Isomorphism Theorem which holds for the
so-called topologically free actions. We present here a definition
of a topological freedom  for complete interactions, which in fact
is a verbatim of the corresponding definition for
partial actions, see \cite{top-free}. Though in the  generality
under consideration we failed to establish the Isomorphism Theorem
we managed to obtain a partial result, see Theorem
\ref{partIsomThm}, and we obtained a complete goal, see Theorem
\ref{complete goal}, in the case of ABL-crossed products, that is
when one of the actions from the pair $(\V,\H)$ acts by
endomorphisms.
\\

The author wishes to express his thanks to  A. V. Lebedev for a number of  comments and remarks which contributed to the preparation of the present paper. 
 \section{Interactions}\label{the_interactions}
Let us start with establishing notation and more accurate
definitions of basic notions appearing in the text. Throughout the
paper we let  $\cal A$ denote a $C^*$-algebra with an identity
$1$, and $\Gamma^+$ to be a positive cone of a totally ordered
abelian group $\Gamma$ with an identity $0$:
$$
\Gamma^+=\{x\in \Gamma: 0\leq x\}, \qquad \Gamma=\Gamma^+\cup( - \Gamma^+), \qquad \Gamma^+\cap( - \Gamma^+)=\{0\}.
$$
\begin{I}
We  say that  $\V$ is an \emph{action} of $\Gamma^+$ on $\A$ if  $\V:\Gamma^+\rightarrow \mathrm{PosLin}(\A)$
is semigroup homomorphisms, and then for each $x\in \Gamma^+$ we denote by $\V_x:\A\to\A$ the corresponding positive linear map:
$$
  \V_0={\rm Id}, \quad \qquad  \V_x\circ \V_y =\V_{x+y}, \qquad \qquad  x,y\in\Gamma^+.
$$
If $\V$ acts not only by linear but also multiplicative maps then
 usually we shall denote it by $\al$ and  call the triple $(\A,\Gamma^+,\al)$ a $C^*$-\emph{dynamical system},
 cf. \cite{Kwa-Leb}.
 \end{I}
 The following is a simple modification of \cite[Definition 3.1]{Kwa-Leb}.
\begin{I}
\label{interactions}
A pair $(\V,\H)$ consisting of two actions $\V$ and $\H$ of $\Gamma^+$
on a $C^*$-algebra $\A$ will be called  \emph{interaction}  if for each $x\in \Gamma^+$
the following conditions are satisfied
  \begin{itemize}
     \item[(i)] $\V_x\H_x\V_x=\V_x$,
  \item[(ii)] $\H_x\V_x\H_x=\H_x$,
  \item[(iii)] $\V_x(ab)=\V_x(a)\V_x(b)$, if either $a$ or $b$ belong to $\H_x(\A)$,
    \item[(iv)] $\H_x(ab)=\H_x(a)\H_x(b)$, if either $a$ or $b$ belong to $\V_x(\A)$.
 \end{itemize}
\end{I}
We  stress that the preceding definition is not a straightforward
generalization of the one given by R. Exel in \cite{exel-inter} (and presented above in the introduction).
\begin{Ee}\label{not power part-iso} Let  $\A=M_2(\C)$ be the algebra of $2\time 2$
complex matrices. We define two positive maps on $\A$ by the formulae
$$
\V((a_{ij}))=\frac{a_{11}}{2}\left(\begin{array}{cc}
1 & 1 \\
1 & 1
\end{array}\right),
\qquad
\H((a_{ij}))=\frac{a_{11}+a_{12}+a_{21}+a_{22}}{2}\left(\begin{array}{cc}
1 & 0 \\
0 & 0
\end{array}\right).
$$
It is a pleasant exercise to show that $\V$ and $\H$ satisfies the conditions
i) - iv) from \ref{interactions}, and hence they form an interaction in the sense of \cite{exel-inter}.
 But they do not yield an interaction in our sense because,
 for instance, $\H^2\circ \V^2 \circ \H^2\neq \H^2$. Actually,
 the obstacle here is that $\V$ and $\H$ are implemented by a partial isometry
 which is not a power partial isometry (in particular $\V(1)\H(1)\neq \H(1)\V(1)$, cf. Proposition \ref{prod-partis}).
\end{Ee}
However,  thanks to \cite[Propositions 2.6 and 2.7]{exel-inter} the following fundamental properties of interactions are true.
 \begin{Pn} \label{TheSubalgebras}
  Let $(\V,\H)$ be an interaction, and let  $x\in \Gamma^+$ be fixed. Then
  \begin{itemize}
  \item[i)] $\V_x(\A)$ and $\H_x(\A)$ are  $C^*$-subalgebras of $\A$,
  \item[ii)] $E_{\V_x}=\V_x\circ \H_x$ is a conditional expectation onto $\V_x(\A)$,
  \item[iii)] $E_{\H_x}=\H_x\circ\V_x$ is a conditional expectation onto $\H_x(\A)$,
  \item[iv)] the mappings
  $$
   \V_x: \H_x(\A) \to \V_x(\A),\qquad  \H_x: \V_x(\A) \to \H_x(\A)
  $$
   are
*-isomorphisms, each being the inverse of the other, and we have
  $
  \V_x = \V_x\circ E_{\H_x}
  $ and $
  \H_x = \H_x\circ  E_{\V_x}.
  $
 \end{itemize}
   \end{Pn}
As the algebra $\A$ considered here is unital we may (for any
interaction $(\V,\H)$) study the elements $\V_x(1)$, $\H_x(1)$,
$x\in \Gamma^+$, which happen to have very useful properties.
   \begin{Pn} \label{Projections}
Let $\A$ contain the unit $1$  and   let $\V$ and $\H$ be actions
on $\A$ forming an interaction $(\V,\H)$.  Then
  \begin{itemize}
  \item[i)] $\{\V_x(1)\}_{x\in \Gamma^+}$ and $\{\H_x(1)\}_{x\in\Gamma^+}$ form  decreasing families of orthogonal projections,
    \item[ii)] for any $a\in \A$ and any $x,$ $y\in \Gamma^+$ such that $y\ge x$  we have
  $$
  \V_y(\H_x(1)a)=\V_y(a\H_x(1))=\V_y(a),\qquad \H_y(\V_x(1)a)=\H_y(a\V_x(1))=\H_y(a),
  $$
  in particular for $y\ge x$
  $$
  \V_y(\H_x(1))=\V_y(1),\qquad \H_y(\V_x(1))=\H_y(1),
  $$
  \item[iii)] for any $a\in \A$ and any $x,$ $y\in \Gamma^+$ such that $y\ge x$  we have
  $$
  a=\V_x(1)a= a\V_x(1), \qquad \textrm{if }\,\,a \in \V_{y}(\A),
  $$
  $$
  a=\H_x(1)a= a\H_x(1), \qquad \textrm{if }\,\,a \in \H_{y}(\A).
  $$
 \end{itemize}
  \end{Pn}

{\bfseries Proof.} ii). Let us observe first that
$E_{\V_x}(1)=\V_x(1)$. Indeed, we have
$$
E_{\V_x}(1)=\V_x(\H_x(1))=\V_x(\H_x(1)1)=\V_x(\H_x(1))\V_x(1)=\V_x(\H_x(1))\V_x\big(\H_x(\V_x(1))\big)
$$
$$
=\V_x\big(\H_x(1)\H_x(\V_x(1))\big)=\V_x\big(\H_x(1\V_x(1))\big)=\V_x\big(\H_x(\V_x(1))\big)=\V_x(1).
$$
For any $a\in\A$ and any $y\ge x$ thus  we have
$$
\V_y(\H_x(1)a)=\V_{y-x}(\V_x(\H_x(1)a))=\V_{y-x}(\V_x(\H_x(1))\V_x(a))=\V_{y-x}(\V_x(1)\V_x(a))=\V_y(a).
$$
Taking adjoints one obtains $\V_y(a\H_x(1))=\V_y(a)$ and hence (by symmetry) ii) is proved.
\\
iii). Let $a =\V_x(b)$ for a certain $b\in \A$. By ii)  we  have
$$
a=\V_x(b)=\V_x(\H_x(1)b)=\V_x(\H_x(1)b)=\V_x(\H_x(1))\V_x(b)=\V_x(1)a,
$$
and similarly  $a= a\V_x(1)$. The proof for $\H_x$ is analogous.
\\
i) Let us  show that $\V_x(1)$ is a projection. Since $\V_x$ is positive
 $\V_x(1)$ is self-adjoint and it is an idempotent because
$$
\V_x(1)=\V_x\big(E_{\H_x}(1)\big)=\V_x\big(E_{\H_x}(1)1 \big)=\V_x\big(E_{\H_x}(1)\big)\V_x(1)=\V_x(1)\V_x(1).
$$
Now let us observe that $\V_x(1)\ge \V_y(1)$ for $y\ge x$. Indeed, using ii) twice we have
$$
\V_x(1)\V_y(1)=\V_x(1)\V_{x}(\V_{y-x}1)=\V_x(\H_x(1))\V_{x}((\V_{y-x}(1))=\V_x(\H_x(1)\V_{y-x}(1))
$$
$$
=\V_x(\V_{y-x}(1))=\V_y(1).
$$
Thus by symmetry $\V_x(1)$ and $\H_x(1)$ form decreasing families of  projections.
\qed
\\

As one would like to think of interactions as of the natural
generalization of $C^*$-dynamical systems, one may be disappointed
to see that for a $C^*$-dynamical system $(\A,\Gamma^+,\al)$ and
its transfer action $L$, the pair $(\al,L)$ might not be an
interaction. However, if the transfer action $L$ is complete the
pair $(\al,L)$ is always an interaction, see Proposition
\ref{complete transfer action},  and the class of transfer actions
that yield  interactions is even wider (for definitions of
transfer and complete transfer actions see Introduction). As an
example we present here a simple corollary to Proposition 3.4 from
\cite{exel-inter}.
\begin{Pn}\label{exel proposition}
Let $L$ be a transfer action for a $C^*$-dynamical system $(\A,\Gamma^+,\al)$
such that $L_x(1)=1$ for each $x\in \Gamma^*$. Then $(\al,L)$ is an interaction.
\end{Pn}

\section{Complete interactions}\label{complette interactions}
Here we introduce a notion of a complete interaction  which is a generalization of the
 complete transfer action  notion, see \ref{complete transfer action}. Afterwards,  for a given  action $\V$ we write down
 the necessary
 and sufficient conditions for existence of an action $\H$  such that $(\V,\H)$ is a complete interaction.
 Moreover, we show the uniqueness of  such action $\H$, see Theorem \ref{complete}.
 In order to show that in general an action does not determine uniquely an interaction
 we adapt to our needs an example from  \cite{Bakht-Leb}.
\begin{Ee}\label{2-adic} (an example of a $C^*$-dynamical system $(\A,\Gamma^+,\al)$ which admits
uncountably many transfer actions  satisfying assumptions of
Proposition \ref{exel proposition}). Let  $\A = C(X)$ where $X =
\R\, (\mathrm{mod}\, 1)$ and let $\Gamma^+=\N$. We define an
action $\al$ by endomorphisms of $\A$   by the formula
\[
\al_n (a) (x) = a (2^nx\, (\mathrm{mod}\, 1)),\qquad n \in \N.
\]
We fix  any continuous function $\rho$ on $X$ having the properties
$$
0\le \rho (x) \le 1, \qquad
\rho \Big(\frac{x}{2}+\frac{1}{2}\Big) +\rho \Big(\frac{x}{2}\Big)= 1, \,\,\,\quad x\in [0,1).
$$
Take  the standard tent map: $T(x)=1-|1-2x|$, $x\in[0,1]$, and
associate with $\rho$ a family of cocycles given by
$$
\rho_0\equiv 1 \,\,  \textrm{ and } \,\, \rho_n(x)=\rho(T^{n-1}(x))\cdot ...\cdot \rho(T(x))\cdot\rho(x), \qquad \textrm{for }n>0.
$$
Then it is not hard to check that $\rho_n$ satisfies the relations
$$
0\le \rho_n (x) \le 1, \qquad
\sum_{k=0}^{2^n-1}\rho_n\Big(\frac{x}{2^n}+\frac{k}{2^n}\Big)= 1,\,\,\, \quad x\in [0,1),
$$
and the following formula defines an action $L$ on $C(X)$
$$
L_n  (a)(x) = \sum_{k=0}^{2^n-1}\rho_n\Big(\frac{x}{2^n}+\frac{k}{2^n}\Big)
 a\Big(\frac{x}{2^n}+\frac{k}{2^n}\Big),\qquad x\in [0,1),\,\,n\in \N.
$$
Clearly for any  $\rho $ chosen $L$ is a transfer action
for $\al$  and since $L_n(1)=1$ for each $n\in \N$, the pair $(\al,L)$ is an
interaction by Proposition \ref{exel proposition}.
\end{Ee}

\begin{Dn}\label{complete transfer action}
The interaction  $(\V,\H)$ will be called {\em complete}, if the following conditions
 are satisfied
\begin{equation}\label{completeness property}
\H_x(\V_x(a))=\H_x(1)a\H_x(1),\quad \V_x(\H_x(a))=\V_x(1)a\V_x(1),\quad \qquad x\in \Gamma^+,\quad a\in \A,
\end{equation}
\begin{equation}\label{completeness property2}
\H_y(1)\V_x(1)=\V_x(1)\H_y(1),\qquad \quad x,y\in \Gamma^+.
\end{equation}
\end{Dn}

The interaction in Example \ref{2-adic} is not complete because  condition
 \eqref{completeness property} is not fulfilled.  The condition  \eqref{completeness property2}
 is closely related to the following  criterium
for the  product of partial isometries to be a partial isometry, cf. Example
 \ref{not power part-iso}, and the proof of Proposition \ref{finely rep}.
\begin{Pn}\label{prod-partis}{\upshape\cite[Lemma 2]{Hal-Wal}}\ \
Let $S$ and
$T$ be partial isometries. Then $ST$ is a partial isometry iff
$S^*S$ commutes with  $TT^*$.
\end{Pn}
Now we explain  the relationship between complete interactions and the complete transfer
actions for $C^*$-dynamical systems.
We denote by $Z(\A)$ the center of $\A$.
\begin{Pn}\label{complete transfer action}
If  $L$  is a  complete transfer action for a $C^*$-dynamical system $(\A,\Gamma^+,\al)$,
 then  the pair $(\al, L)$ is a complete interaction
 and
$$
L_x(1)\in Z(\A),\qquad  x\in \Gamma^{+}.$$
Conversly, if $(\V,\H)$ is a complete interaction such that $\H_x(1)\in Z(\A)$,
$x\in \Gamma^{+}$, then $(\A,\Gamma^+,\V)$ is a $C^*$-dynamical system and $\H$
is its complete transfer action.
\end{Pn}
{\bfseries Proof.} Let us prove the first part of the proposition.
\ref{interactions}.i) follows from \cite[2.2]{Kwa-Leb}, and  \ref{interactions}.ii)
follows from \cite[2.3]{Kwa-Leb}, see also \cite[(2.15)]{Bakht-Leb}. Since $\al_x$
is an endomorphism \ref{interactions}.iii) is trivial.  We recall that $L_x(1)$ belongs
to the center  of $\A$ and $L_x(\al_x(a))=L_x(1)a$, cf. \cite[Theorem 2.4]{Kwa-Leb}.
Hence \eqref{completeness property}, \eqref{completeness property2} are valid and to show
\ref{interactions}.iv)  we notice that
$$
L_x(\al_x(a)b)=aL_x(b)=aL_x(1)L_x(b)=L_x(\al_x(a))L_x(b).
$$
By taking adjoints one obtains $L_x(b\al_x(a))=L_x(b)L_x(\al_x(a))$.
\\
To prove the remaining part of the statement it suffices to show that if $\H_x(1)$
belongs to $Z(\A)$ then $\V_x$ is multiplicative. By
Proposition \ref{Projections}, formula \eqref{completeness property}
and the definition of interaction  we  have
$$
\V_x(ab)=\V_x(\H(1)ab\H(1))=\V_x(a\H(1)b(1)\H(1))=\V_x(a)\V_x(\H(1)b(1)\H(1))=\V_x(a)\V_x(b)
$$
for arbitrary $a,b\in A$, and the proof is complete.
\qed
\\
In view of the above proposition the following statement is a
generalization of  \cite[Theorem 2.4]{Kwa-Leb}.
\begin{Tm}\label{complete}
Let $\V$ be an action of $\Gamma^+$ on $\A$. The following are equivalent:
\begin{itemize}
\item[$1)$]   there exists an action  $\H$   such that $(\V,\H)$ is a complete interaction,
\item[$2)$] $(i)$ there exists an action  $\H$   such that $(\V,\H)$ is an interaction,\\[6pt]
$(ii)$   $\V_x(\A)$, $\H_x(\A)$  are   hereditary   subalgebras of $\A$ for each $x\in\Gamma^+$,\\[6pt]
$(iii)$ $\V_x(1)$ and $\H_y(1)$ commute for all $x,y\in\Gamma^+$,
\item[$3)$] $(i)$  $\V_x(1)$ is an orthogonal projection and $\V_x(\A)=\V_x(1)\A\V_x(1)$ for each $x\in\Gamma^+$,\\[6pt]
$(ii)$ there exists a decreasing family $\{P_x\}_{x\in\Gamma^+}$ of  orthogonal projections such that
\begin{itemize}
\item[$a)$] $\V_x (1)$ and $P_{y}$  commute for all  $x,y\in\Gamma^+$,
\item[$b)$] $\V_x(P_{x+y}) = \V_x (1)P_y$, for each $x,y\in\Gamma^+$,
\item[$c)$] the mappings $\V_x \!: P_x \A P_x \to
\V_x (\A)$ are $^*$-isomorphisms.
\end{itemize}
\end{itemize}
Moreover the objects in $1)$ -- $3)$ are defined in a unique way,
i.e. the  action $\H$ in $1)$ and $2)$ is
unique and the family  of projections $\{P_x\}_{x\in\Gamma^+}$ in $3)$ is unique as well. These object are combined by formulae
\begin{equation}
\label{P} P_x=\H_x(1),\qquad x\in\Gamma^+,
\end{equation}
and
\begin{equation}\label{d*-}
\H_x(a) =\V_x^{-1}(\V_x(1)a\V_x(1)), \ \ \  a\in \A, \ \
\end{equation}
where $\V_x^{-1}:\al_x(\A)\to P_x\A P_x$ is the inverse
mapping to \ $\V_x:P_x\A P_x \to \V_x(\A)$,  $x\in \Gamma^+$.
\end{Tm}
{\bfseries Proof.} $1)\Leftrightarrow 2)$. In view of
\eqref{completeness property} and Proposition \ref{TheSubalgebras}
it is enough to show that  2) (ii) is equivalent to
\eqref{completeness property2}. It is straightforward that if
\eqref{completeness property2}  holds, then
$$
\H_x(\A)=\H_x(1)\A \H_x(1),\qquad \V_x(\A)=\V_x(1)\A\V_x(1)
$$
are hereditary subalgebras. Conversely, if $\H_x(A)$ and $\V_x(\A)$ are hereditary subalgebras
of $\A$, then the argument used in the proof of  \cite[Proposition 4.1]{exel} shows
that $\V_x(1)\A \V_x(1)\subset \V_x(\A)$  and $\H_x(1)\A \H_x(1)\subset \H_x(\A)$.
By Proposition \ref{Projections} we have  $\H_x(\A)\subset H_x(1)\A \H_x(1)$ and
$\V_x(\A)\subset\V_x(1)\A\V_x(1)$, and hence  \eqref{completeness property2} holds.
\\
$1), 2)\Rightarrow 3)$. Take $P_x=\H_x(1)$, $x\in \Gamma^+$. Item 3) then follows from
 Propositions \ref{TheSubalgebras} and  \ref{Projections}.
\\
$3) \Rightarrow 1)$. Fix $x\in \Gamma^+$. Let
$\V_x^{-1}\!:\V_x(\A)=\V_x(1)\A\V_x(1)\to P_x\cal A P_x$ be the
inverse mapping to  $\V_x\!:P_x\A P_x\to \V_x(\A)$. Define
$\H_x$ by the formula $\H_x(a) =\V_x^{-1}(\V_x(1)a\V_x(1))$.
Clearly $\H_x$ is linear and  positive, and (\ref{completeness
property}) is fulfilled. Furthermore, \ref{interactions}.i), ii)
hold. To prove \ref{interactions}.iii) we  note that
$$
\V_x\big(\H_x(\V_x(a)b)\big) =\V_x(1)\V_x(a)b\V_x(1) =\V_x(a)\V_x(1) b\V_x(1) =\V_x(a)\V_x(\H_x(b))
$$
$$
=\V_x(\H_x(\V_x(a)))\V_x(\H_x(b))=\V_x( \H_x(\V_x(a)) \H_x(b)),
$$
and as the elements  $\H_x(\V_x(a)b)$
and
$\H_x(\V_x(a))\H_x(b)$ belong to the subalgebra  $P_x\A P_x$ where the mapping $\V_x$
is injective,  they coincide. Similarly one proves that $\H_x(a\V_x(b))=\H_x(a)\H_x(\V_x(b))$
and thus \ref{interactions}.iii) holds.
\\
The same argument proves \ref{interactions}.iv) and
therefore to show that $(\V,\H)$ is an interaction we only need to prove  that
 $\H$ is an action of the semigroup $\Gamma^+$.
\\
Using $3)$ $(ii)$  and  \ref{interactions}.iii)  we have
$$
\V_y(P_{x+y}\A P_{x+y})=\V_y(P_{x+y})\V_{y}(\A)\V_y(P_{x+y})= P_x\V_y(\A)P_x
$$   and as $P_{x+y}\A P_{x+y}\subset P_y\A P_{y}$ we obtain  that
$\V_y:P_{x+y}\A P_{x+y}\rightarrow P_x\V_y(\A)P_{x}$ is a $^*$-isomorphism and the inverse is given by $\H_y$.
 Thus    we have
$$
\H_y(\H_x(\A))=\H_y(P_x\A P_x)=\H_y(\V_y(1)P_x\A P_x\V_y(1))
$$
$$
=\H_y(P_x\V_y(1)\A\V_y(1) P_x)=\H_y(P_x\V_y(\A) P_x)=P_{x+y}\A P_{x+y}.
$$
Hence $\H_y(\H_x(a))$ and $\H_{x+y}(a)$ belong to the subalgebra  $P_{x+y}\A P_{x+y}$ where the map $\V_{x+y}$
is injective, and  as
$$
\V_{x+y}(\H_y(\H_x(a))=\V_{x}\big(\V_y(\H_y(\H_x(a))\big)=\V_{x}\big(\V_y(1)\H_x(a)\V_y(1)\big)
$$
$$
=\V_{x}\big(\V_y(1)P_x\H_x(a)P_x\V_y(1)\big)=\V_{x}\big(P_x\V_y(1)P_x\H_x(a)P_x\V_y(1)P_x\big)
$$
$$
=\V_{x}(P_x\V_y(1)P_x)\V_x(\H_x(a))\V_x(P_x\V_y(1)P_x)
= \V_{x}(\V_y(1))\V_x(\H_x(a))\V_{x}(\V_y(1))
$$
$$
= \V_{x+y}(1)\V_x(1)a\V_x(1)\V_{x+y}(1)
=\V_{x+y}(1)a\V_{x+y}(1)=\V_{x+y}(\H_{x+y}(a))
$$
we have  $\V_{x+y}=\V_y\circ \V_x$.
\\
The uniqueness of the objects in 1) - 3) is  straightforward. \qed

\section{Finely representable actions and their \\crossed products }
\label{endom}
In this section  we define  finely representable actions as the ones possessing nondegenerated covariant representations, and thereby
 possessing nondegenerated crossed products.  These actions are closely related to complete interactions.
  Namely, it is not very difficult to prove (see Proposition \ref{finely rep}) that every finely
  representable action is a 'part' of a complete interaction, and although it might be  difficult
  to prove it is very likely that the opposite is also true, cf. \cite{Bakht-Leb}, \cite{Kwa-Leb}.
\\
Furthermore, we investigate a dense $^*$-subalgebra of the crossed product via  quasi-mono\-mials.
In particular we prove certain inequality which will be of  primary importance  in the forthcoming sections.

\begin{Dn}
\label{fine}
Let  $\V$ be  an action of $\Gamma^+$ on a $C^*$-algebra  $\A$. We say that  $\V$ is
{\em finely representable\/} if
there exists a triple $(C,\sigma ,U)$, called a \emph{covariant representation} of $\V$,
consisting of a unital $C^*$-algebra $C$, unital monomorphism  $\sigma:\A\to
C$ and a semigroup homomorphism $U:\Gamma^+ \to C$   such that
 for every $x\in \Gamma^+$, $U_x$ is a partial isometry, and for every  $a\in\A$, $x\in\Gamma^+$, the
following conditions are satisfied
\begin{equation}\label{b,,4}
\sigma(\V_x(a)) =U_x\sigma(a)U^*_x,\qquad U^*_x\sigma(a)U_x \in \sigma(\A)
\end{equation}
\end{Dn}
Let us clarify how  the interaction notion is involved in the
above definition.
\begin{Pn}\label{finely rep}
If $\V$ is a finely  representable action of $\Gamma^+$ on  $\A$,
then there exists a (necessarily unique) action $\H$ such that
$(\V,\H)$ is a complete interaction. Moreover for any covariant
representation $(C,\sigma ,U)$ the following formulae hold
\begin{equation}\label{takie tam equation}
\sigma(\V_x(a))=U_x\sigma(a)U_x^*,\qquad \sigma(\H_x(a))=U_x^*\sigma(a)U_x, \qquad \quad a\in \A, \, x\in \Gamma^{+}.
\end{equation}
\end{Pn}
{\bfseries Proof.}
If conditions (\ref{b,,4})
  are satisfied then (identifying $\cal A$ with
$\sigma (\cal A)$) one can set
\[
\H_x (\cdot ) = U_x^* (\cdot )U_x, \qquad x\in \Gamma^+.
\]
Using fundamental properties of partial isometries one easily
verifies that $(\V,\H)$ is an interaction and that conditions
\eqref{completeness property} are  satisfied. Condition
\eqref{completeness property2} follows from the fact that $U_x U_y=
U_{x+y}$ is a partial isometry, and Proposition \ref{prod-partis}.
Thus  $(\V,\H)$ is a complete interaction. By Theorem
\ref{complete},  $\H$ is unique, and hence \eqref{takie tam
equation} holds for any covariant representation of $\V$. \qed
\\
The following statement  is partially converse to the above one.

\begin{Tm}\label{C-dynamical conditions}
Let $(\V,\H)$ be a complete interaction such that one of the
equivalent conditions $i)$, $ii)$, $iii)$ hold
\begin{itemize}
\item[i)] each $\V_x$ is an endomorphism,
\item[ii)] $\H_x(1)\in Z(\A)$, for all $x\in \Gamma^+$,
\item[iii)] $(\A,\Gamma^+,\V)$ is a $C^*$-dynamical system,
\end{itemize}
or a counter part of one of them  with $\V$  replaced by $\H$ hold.
Then both $\V$ and $\H$ are finely representable actions.
\end{Tm}
{\bfseries Proof.}
If follows from Proposition \ref{complete transfer action} and \cite[Theorem 3.2]{Kwa-Leb}.
\qed
\par
Unfortunately the author was not able to answer the following general question:
\begin{quote}
\textbf{Problem. } Let $(\V, \H)$ be an arbitrary complete interaction. Are the actions $\V$ and $\H$ finely representable?
\end{quote}
Fortunately, this obstacle does not really affect our further considerations.
\\
Let us note  that by Proposition \ref{finely rep} every finely representable action $\V$
 determines uniquely another finely representable action $\H$  such that for every
 covariant representation $(C,\sigma,U)$  of $\V$ the triple $(C,\sigma,U^*)$ where $(U^*)_x=U_x^*$,
  is a covariant representation for $\H$ and vice versa. In particular, $\H$ is finely representable
  and in view of the following definition the crossed products by $\V$ and $\H$ coincide.

 \begin{Dn}\label{crossed-product definition}
Let $\V$ be a finely representable action and let $(\V,\H)$ be the
corresponding complete interaction. The {\em crossed product}
(also called  \emph{covariance algebra}) of  the $C^*$-algebra
$\A$ by the action $\V$, which we denote by $\A\times_{(\V,\H)}
\Gamma$ to indicate the role (and the symmetry) of the interaction
$(\V,\H)$, is the universal unital $C^*$-algebra generated by a
copy of $\A$ and a family $\{\hat{U}_x\}_{x\in \Gamma^{+}}$ of
partial isometries  subject to relations
\begin{equation}
\label{,b,,4}
\V_x(a) = \hat{U}_xa\hat{U}_x^*,\quad\ \H_x (a)=\hat{U}_x^*a\hat{U}_x ,\qquad a\in \cal A,\,\, x\in \Gamma^{+},
\end{equation}
$$
\hat{U}_x\hat{U}_y=\hat{U}_{x+y}, \qquad x,y \in \Gamma^{+}.
$$
 If $(C,\sigma ,U)$ is a covariant representation of $\V$  then we
 denote by $(\sigma \times U)$ the homomorphism of $\A\times_{(\V,\H)}\Gamma$ into $C$ established by
$$
(\sigma \times U)(a)=\sigma(a), \quad \qquad (\sigma \times U)(\hat{U}_x)=U_x,\quad \quad a\in\A, \,\,x\in \Gamma^{+}.
$$
\end{Dn}

In order to  study   covariance algebras it is important to understand the structure of
 a  $^*$-subalgebra $C_0$ of $\A\times_{(\V,\H)}\Gamma$ generated by $\A$ and a semigroup
 $ \hat{U}=\{\hat{U}_x\}_{x\in \Gamma^+}$. Let us thus investigate $C_0$.
\\
The basic elements in $C_0$ are the ones of the form
\begin{equation}\label{monomials}
\prod_{i=1}^n a_i\hat{U}^*_{x_i}=a_1\hat{U}^*_{x_1}a_2\dots a_n\hat{U}^*_{x_n},\qquad
 \qquad  \prod_{i=1}^n a_i\hat{U}_{x_i}=a_1\hat{U}_{x_1}a_2\dots a_n\hat{U}_{x_n},
\end{equation}
$x_1,...,x_n\in \Gamma^+$,  $a_{1},...,a_{n}\in \A$. We shall call them \emph{monomials}
 of negative  and positive type respectively. In this context the element  $x_1+...+x_n$
 is  a \emph{degree} of both of  these  monomials, and any  finite sum of monomials of the
 same type and the same degree will be called a \emph{quasi-monomial}.
 Namely quasi-monomials of degree $x$ are the elements of the form
\begin{equation}\label{quasi-monomials}
q_{-x}=\sum_{y=(y_1,...,y_n)\in Q \atop  y_1+...+y_n=x}\prod_{i=1}^n a_i^{-y}\hat{U}^*_{y_i},
\qquad \qquad q_x=\sum_{y=(y_1,...,y_n)\in Q\atop y_1+...+y_n=x }\prod_{i=1}^n a_i^{y}\hat{U}_{y_i}
\end{equation}
where $Q$ is a finite set consisting of finite  sequences with
entries in $\Gamma^+$ (presumably with different lengths). In
particular every quasi-monomial $q_0$  of degree $0$ is in fact a
monomial and $q_0\in \A$.
 We claim that
\begin{Pn}
  $C_0$ consists of  finite sums of  monomials \eqref{monomials},
  and a fortiori  of sums of quasi-monomials.
\end{Pn}
{\bfseries Proof.} It is clear that the finite sums of monomials
form  a self-adjoint linear space (containing $\A$ and
$\{\hat{U}_x\}_{x\in\Gamma}$). In fact they form an algebra
because every "mixed monomial"
$a_1\hat{U}_{x_1}b_1\hat{U}_{y_1}^{*}a_2\hat{U}_{x_2}\dots
a_n\hat{U}_{x_n}b_n{\hat{U}_{y_n}^{*}}$ equals to  a "non-mixed
monomial" in one  of the forms
$$
c_1\hat{U}^*_{z_1}c_2\dots c_m\hat{U}^*_{z_m}\qquad
 \textrm{ or }\qquad  c_1\hat{U}_{z_1}c_2\dots c_m\hat{U}_{z_m}
$$
depending on whether $x_1+\dots + x_n\leq y_1+\dots + y_n$ or
$y_1+\dots + y_n \leq x_1+\dots + x_n$ (this is an easy fact due to the
total ordering of $\Gamma$ and  \eqref{,b,,4}.
\qed
\\
Consequently, for any $a\in C_0$ there exists a finite set
$F\subset \Gamma^{+}\setminus \{0\}$, and a family of
quasi-monomials $q_{\pm x}$ of degree $x\in F$and $a_0\in \A$,
such that
\begin{equation}\label{elements of the form}
a=\sum_{x\in F} q_{-x} + a_0 +\sum_{x \in F} q_x.
\end{equation}
Moreover, as the next proposition shows, the quasi-monomial $a_0$ of degree $0$ is uniquely determined by $a$.
\begin{Pn}\label{B0}
For any $a\in C_0$, and any presentation of $a$ in the form \eqref{elements of the form}
the following inequality holds
\begin{equation}\label{star property}
\|a_0\|\leq \|a\|.
\end{equation}
\end{Pn}
\textbf{Proof.}
Take any   faithful non-degenerate representation  $\pi:\A\times_{(\V,\H)}\Gamma \to H$,
 and consider the Hilbert space $\widetilde{ H} =\bigoplus_{g\in \Gamma} H_g$ where $H_g=H$,
 for all $g\in \Gamma$, and the
representation  $\nu : \A\times_{(\V,\H)}\Gamma\to L(\widetilde{
H})$  given by the formulae
\begin{gather*}
(\nu (a)\xi )_g = \pi (a) (\xi_g), \qquad\text{where}\quad
a\in {\cal A}, \quad \widetilde{ H} \ni \xi = \{ \xi_g  \}_{g\in
\Gamma}\,;\\[6pt]
(\nu (\widetilde{U}_x)\xi )_g =\pi(\widetilde{U}_x) (\xi_{g-x}),\qquad
(\nu (\widetilde{U}^*_x)\xi )_g= \pi (\widetilde{U}^*_x) (\xi_{g+x}).
\end{gather*}
Routine verification shows that $\nu (\A)$ and $\nu
(\widetilde{U}_x)$ satisfy all the conditions of a covariant representation and thus
 $\nu$ is well defined.
\\
Now take any $a\in  \A\times_{(\V,\H)}\Gamma$ given by
\eqref{elements of the form} and for a given $\varepsilon > 0$ chose  a vector
$\eta \in H$ such that
\begin{equation}\label{e}
\Vert \eta \Vert =1 \quad \text{and}\quad \Vert \pi (a_0)
\eta \Vert > \Vert \pi (a_0)  \Vert - \varepsilon .
\end{equation}
Set $\xi \in  \widetilde{ H}$ by $\xi_g = \delta_{(0,g)}\eta $,
where $\delta_{(i,j)}$ is the Kronecker symbol. Then we have  $\Vert
\xi \Vert = 1$ and the explicit form of $\nu (a) \xi$  and
\eqref{e} imply
\[
\Vert \nu (a) \xi \Vert \ge \Vert \pi (a_0) \eta \Vert >
\Vert \pi (a_0)  \Vert - \varepsilon
\]
which by the arbitrariness of $\varepsilon$ proves the desired
inequality:
$$
\Vert a  \Vert \ge \Vert \nu(a)  \Vert\ge \Vert \pi (a_0)  \Vert = \Vert  a_0 \Vert.
$$
\qed
\begin{Rk}\label{borek jest wolny}
It is clear that the form  \eqref{elements of the form} of $a\in C_0$
is far from being  unique in general. However, if $(\V,\H)$ comes from a
$C^*$-dynamical system, i.e.  one of the conditions $i)$-$iii)$ from
Theorem \ref{C-dynamical conditions}  holds, then
 every monomial and every quasi-monomial of degree $x\in \Gamma^{+}$
can be presented in one of the  forms $q_{-x}=\hat{U}_x^*a_{-x}$ or $q_x=a_{x}\hat{U}_x$,
 cf. \cite{Leb-Odz}. Consequently, any element $a\in C_0$  can be presented in the form
$$
a=\sum_{x\in F} \hat{U}_x^*a_{-x}  + a_0 + \sum_{x\in F} a_{x}\hat{U}_x
\quad
\textrm{ where }\,\,a_{-x}\in  \A \hat{U} \hat{U}^* \,\,\textrm{  and } \,\, a_{x}\in   \hat{U}\hat{U}^*\A .$$
Moreover, see \cite{Kwa-Leb}, \cite{Leb-Odz}, the coefficients $a_{\pm x}$ in the above formula are uniquely determined by $a$.
\end{Rk}

\section{Conditional expectation and
\\faithful representations of  crossed products}
\label{2.1}
From now on, we fix a finely representable action $\V$ and hence by Proposition
\ref{finely rep} we also fix a complete interaction $(\V,\H)$.
Here we use Proposition \ref{B0} to define a conditional expectation from
$\A \rtimes_{(\V,\H)} \Gamma$ onto $\A$, for which certain 'spectral' formula holds,
 see \eqref{be3.131},  and to give a criterion  for a representation of
 $\A \rtimes_{(\V,\H)} \Gamma$ to be faithful. In the literature
 such necessary and sufficient condition plays important role and is usually called
  property $(^*)$ (for different versions and a history of property $(^*)$ see in particular \cite{Leb-Odz},
  \cite{top-free}, \cite{Ant-Bakht-Leb}, \cite{kwa}).
 \par
The first advantage  of inequality \eqref{star property} is that it implies that the mapping $E_0:C_0\to \A$ given by
$$
 E_0(a)= a_0
$$
where $a$ is of the form \eqref{elements of the form}, is well defined and can be
extended to the conditional expectation acting on the whole of $\A
\rtimes_{(\V,\H)} \Gamma$. We shall show that using $E_0$  one may
express (by the formula generalizing the $C^*$-equality
$\|a\|^2=\|aa^*\|$, $a\in \A$) the norm of  elements from $\A
\rtimes_{(\V,\H)} \Gamma$ by the norms of elements from $\A$, see
Theorem \ref{3a.N}. But first we need to estimate the growth rate
of number of quasi-monomials  appearing in the powers of an
element $a\in C_0$.
\begin{Pn}\label{growth rate}
For any  $a\in C_0$ there exists a family $\{F_k\}_{k\in \N}$ of finite
 subsets of $\Gamma^{+}\setminus\{0\}$ such that
$$
a^k=\sum_{x\in F_k} q_{-x}(k) + q_0(k) +\sum_{x \in F_k} q_x(k)
$$
where $q_{\pm x}(k)$ are quasi-monomials of degree $x$, $x\in \Gamma^+$, $k\in \N$, and
$$
\lim_{k\to \infty} |F_k|^{\frac{1}{k}}=1
$$
where $|F|$ denotes the number of elements in a set $F$. In other
words, the growth rate of number  of quasi-monomials
appearing in the $k$-th power of $a$ is subexponential.
\end{Pn}
{\bfseries Proof.} Let $a$ be given by \eqref{elements of the
form} where $F=\{x_1,..,x_n\}$, then the quasi-monomials in
\eqref{elements of the form} are numbered by the elements of
$F_0=\{0, \pm x_1, ...,\pm x_n\}$, and it is clear that the
quasi-monomials appearing in $a^k$ may be numbered by the set
$F_0^k=\{y_1y_2...y_k: y_i\in F_0\}$. Thus putting $F_k=F_0^k\cap (\Gamma^{+}\setminus\{0\})$ and
recalling that  abelian groups are  subexponential   one obtains the
hypotheses. \qed
\begin{Tm}
\label{3a.N} Let $a\in C_0\subset \A \rtimes_{(\V,\H)} \Gamma$. Then we have
\begin{equation}
\label{be3.131} \Vert a \Vert = \lim_{k\to\infty} \sqrt[4k]{ \left\Vert E_0 \left[
(a\cdot a^*)^{2k}\right]\right\Vert }.
\end{equation}
\end{Tm}
{\bfseries Proof.}  Let $a$ be of the form \eqref{elements of the form}.
 Applying to $a$ the known equality $\left\Vert  \sum_{i=1}^m d_i  \right\Vert^2 \le m \left\Vert  \sum_{i=1}^m d_i
d_i^* \right\Vert$ (which holds for any elements $d_1, ..., d_m$ in an arbitrary $C^*$-algebra) where  $m=2|F|+1$ and $
d_k$, $k=1,...,m$,  are appropriate quasi-monomials,
 we obtain that
 $$\Vert a \Vert^2 \leq
(2|F|+1) \left\Vert  a_0a_0^* + \sum_{x\in F} \big(q_{-x}q_{-x}^* + q_{x} q_{x}\big)\right\Vert= (2|F|+1) \Vert E_0 (aa^*)\Vert.
$$
  On the other hand as
$E_0$ is contractive we have
$
 \Vert a \Vert^2 = \Vert aa^* \Vert \ge \Vert E_0
(aa^*)\Vert
$ and
 thus
\begin{equation}
\label{be3.101} \Vert E_0 (aa^*)\Vert \le  \Vert aa^* \Vert = \Vert a \Vert^2  \le (2|F|+1) \Vert
E_0 (aa^*)\Vert
\end{equation}
Applying (\ref{be3.101}) to $(aa^* )^k$ and having in mind that $(aa^* )^k = (aa^* )^{k*} $ and
$\Vert (aa^* )^{2k}\Vert = \Vert a \Vert^{4k} $ one has
$$
\Vert E_0 \left[
(aa^*)^{2k}\right]\Vert \le
 \Vert (aa^*)^k \cdot  (aa^*)^{k*}  \Vert =
\Vert a \Vert^{4k}  \le (2|F_k|+1) \Vert E_0 \left[ (aa^*)^{2k}\right]\Vert
$$
where $F_k \subset\Gamma^+\setminus\{0\}$ is the  set of all
degrees of non-zero quasi-monomials appearing in  $a^k$. By
Proposition \ref{growth rate} we  have  $\lim_{k\to\infty}
(2|F_k|+1)^{\frac{1}{k}}=1$, and thus
 $$
 \sqrt[4k]{ \left\Vert E_0 \left[ (aa^*)^{2k}\right]\right\Vert }\le \Vert a \Vert \le
\sqrt[4k]{ 2|F_k|+1 }\cdot \sqrt[4k]{  \left\Vert E_0 \left[ (aa^*)^{2k}\right]\right\Vert }
 $$
implies that $
   \Vert a
\Vert = \lim_{k\to\infty} \sqrt[4k]{ \left\Vert E_0 \left[ (aa^*)^{2k}\right]\right\Vert }.
 $
\qed
\\
 One would perceive the origin of the following definition in Proposition \ref{B0}.
\begin{Dn} Let $(C,\sigma ,U)$ be a covariant
representation of $\V$. We shall say that $(C,\sigma ,U)$ possesses
{\em property } $(^*)$ if  for any element $a\in C_0$ (that is for $a$ of the form  \eqref{elements of the form})
we have
$$
\|E_0(a)\| \leq \|(\sigma\times U)(a)\|\qquad (^*)
$$
or in other words
$$
\Vert a _0\Vert \leq \Vert \sum_{x\in F} \left(\sum_{y=(y_1,...,y_n)\in Q \atop
 y_1+...+y_n=x}\prod_{i=1}^{n}\sigma(a^{-y}_i)U_{y_i}^*\right) + \sigma(a_0) + \sum_{x\in F}
 \left(\sum_{y=(y_1,...,y_n)\in Q \atop  y_1+...+y_n=x}\prod_{i=1}^{n}\sigma(a^{-y}_i)U_{y_i}\right)\Vert.
$$
\end{Dn}
We are ready to formulate and prove the main result of this section.
\begin{Tm}\label{isomorphism theorem}
Let $(C,\sigma ,U)$ be a covariant representation of $(\V,\H)$.
The homomorphism  $(\sigma\times U):\A \rtimes_{(\V,\H)} \Gamma
\to C$ is a monomorphism if and only if  $(C,\sigma ,U)$ possesses
property $(^*)$.
\end{Tm}
{\bfseries Proof.} Necessity follows from Proposition \ref{B0}. In
order to  show the sufficiency take any $a\in C_0$. By Theorem
\ref{3a.N} and the definition of property $(^*)$ we have
$$
\label{be3.013} \Vert a \Vert = \lim_{k\to\infty} \sqrt[4k]{ \left\Vert E_0 \left[
(aa^*)^{2k}\right]\right\Vert }\leq \lim_{k\to\infty} \sqrt[4k]{ \left\Vert (\sigma\times U)
(aa^*)^{2k}\right\Vert }
$$
$$
=\lim_{k\to\infty} \sqrt[4k]{ \left\Vert (\sigma\times U)
(aa^*)^{k}(\sigma\times U)
(aa^*)^{k}\right\Vert }=\lim_{k\to\infty} \sqrt[4k]{ \left\Vert (\sigma\times U)
(a)\right\Vert^{4k} }=\left\Vert (\sigma\times U)
(a)\right\Vert.
$$
Hence $
\|a \| =\|(\sigma\times U) (a) \| $ on a dense subset of  $\A \rtimes_{(\V,\H)} \Gamma$.
\qed\begin{Cy}\label{dual action}
There is  the action of the dual group $\hat{\Gamma}$ by the automorphisms of $\A \rtimes_{(\V,\H)} \Gamma$   given by
$$
\lambda a:=a, \quad a\in \A, \qquad  \lambda \hat{U}_x :=\lambda_x \hat{U}_x,\qquad
$$
for $x\in \Gamma^{+},\, \lambda\in \hat{\Gamma},\, \lambda_x=\lambda(x)$
(here we consider $\Gamma$ as a discrete group).
\end{Cy}
{\bfseries Proof.} Suppose that  $\A \rtimes_{(\V,\H)} \Gamma$ is
faithfully and nondegenerately represented on a Hilbert space $H$.
Then for each $\lambda \in\hat{ \Gamma}$ the triple $(id, \lambda
\hat{U} , H)$, where $\lambda \hat{U}=\{ \lambda \hat{U}_x\}_{x\in
\Gamma^{+}}$, is a covariant representation possessing property
$(^*)$, whence $(id\times \lambda \hat{U})$ is an automorphism of
$\A \rtimes_{(\V,\H)} \Gamma$. \qed

\section{Topologically free interactions}\label{topolo free}
In this section we rely heavily on the paper \cite{top-free} where A. V. Lebedev undertook
the issue
of topological freedom for partial actions of groups and obtained the Isomorphism Theorem for partial crossed products.
Roughly speaking, the  contribution of the author of the present paper to the current
 section reduces nearly only to an observation that the definition of topological freedom given by
 A. V. Lebedev  also makes sense in the context of complete interactions. In particular,
 Lemma \ref{2.6} and its proof is an almost faithful verbatim  of  \cite[Lemma 2.7]{top-free}.
\par
To start
with let us  note that  a complete interaction defines in a natural way
partial dynamical systems  (the actions of a group $\Gamma$ by
partial homeomorphisms) on the primitive ideal space ${\rm Prim}\,
\A$ and the spectrum $\hat{\A}$ of $\A$ considered here as
topological spaces equipped with  the Jacobson topology.
\\
 Let us give  the description of these partial dynamical systems. For any $x\in \Gamma^+$ we set
$$
\A_x = \V_x(1) \A \V_x(1) ,\qquad \A_{-x}=\H_x(1) \A \H_x(1),
$$
and  thus  we have a family $\{\A_g\}_{g\in \Gamma}$ of hereditary subalgebras of $\A$.
\\
We recall that for any subset  $S\subset \A$ the  set   $\supp S = \{ I\in \Prim \A : x
  \not\supset  S \}$  is  open  in $\Prim  \A$ (see \cite[Proposition 3.1.2]{Dixmier}),
  and for any hereditary $C^*$-subalgebra $\B$ of $\A$
the  mapping  $I \to I\cap \B$ establishes a homeomorphism  $\supp \B \longleftrightarrow
\Prim \B$ (see \cite[Theorem 5.5.5]{Murphybook}). Analogously, the  set $\hat{\A}^S =
 \{ \pi \in \hat{\A} : \pi (S) \neq 0  \}$ is open in $\hat{\A}$ and for any hereditary
 $C^*$-subalgebra $\B$ of $\A$ the mapping  $\pi \to \left. \pi  \right|_\B$ establishes a homeomorphism
$\hat{\A}^\B\longleftrightarrow \hat{\B}$  (see  \cite{Dixmier}, 3.2.1.).
Thus we may  and we shall identify the family $\{\Prim \A_g\}_{g\in \Gamma}$ with
the family  $\{\supp \A_g\}_{g\in \Gamma}$ of open sets in $\Prim \A$, and family
$\{\hat{\A}_g\}_{g\in \Gamma}$ with the family $\{\hat{\A}^{\A_g}\}_{g\in \Gamma}$ of open sets in $\hat{\A}$.
\\
Let us define the mappings  $\tau_x : \hat{\A}_{-x} \to \hat{\A}_{x}$, $\tau_{-x} :
\hat{\A}_{x} \to \hat{\A}_{-x}$, for $x\in \Gamma^+$, in the following way:
$$
\tau_x (\pi )(a)  = \pi (\H_x(a)), \ \qquad \   \pi \in \hat{\A}_{-x},\  \ a\in \A_x,
$$
$$
\tau_{-x} (\pi )(a)  = \pi (\V_x(a)), \ \qquad \   \pi \in \hat{\A}_{x},\  \ a\in \A_{-x}.
$$
By Theorem 5.5.7 in \cite{Murphybook},    $\tau_x$ and $\tau_{-x}$ are homeomorphisms.
Let us also define the mapping $t_g : \Prim \A_{-g} \to \Prim \A_{g}$,  for $g\in \Gamma$  in the following way:
for any point  $I \in \Prim \A_{-g}$ such that  $I = \ker \pi$ where $\pi \in  \hat{\A}_{-g}$
we set
$$
t_g (I) = \ker  \tau_g (\pi ).
$$
Clearly \/$t_g$\/ is a homeomorphism.
\\
Concluding, for $\tau_g$ and $t_g$ defined in the above described
way  $\{ \tau_g  \}_{g\in \Gamma}$ defines an action of $\Gamma$
by  {\em partial homeomorphisms} of $\hat{\A}$ and $\{  t_g
\}_{g\in \Gamma}$  defines an action of $\Gamma$  by  {\em partial
homeomorphisms} of $\Prim \A$.
\begin{Dn}
\label{top-free} We say that the interaction $(\V,\H)$  is {\em
topologically free} if one of the following equivalent conditions
holds
\begin{itemize}
\item[i)] for any finite set
$\{ x_1 , ..., x_k \} \subset \Gamma^{+}$ and any nonempty open set
  $U \subset \Prim \A_{-x_1}\cap ... \cap   \Prim \A_{-x_k}$
there exists a point $I \in U $ such that all the points $t_{x_i}(I)$, $i=1,...,k$ are distinct.

\item[ii)] for any finite set
$\{ x_1 , ..., x_k \} \subset \Gamma^+$ and any nonempty open set $U$ there exists
a point $I \in U $ such that all the points $t_{x_i}(I)$,  $i=1,...,k$  that
 are defined   are distinct.
\item[iii)]
If we denote by $G_x$ the set
$$
G_x = \{I \in \Prim \A_{-x}: t_x (I)=I   \}
$$
for any finite set $\{ x_1, ... x_n \}\subset \Gamma^{+}\setminus
\{0\}$,  the interior of the set $
  \bigcup_{i=1}^n G_{x_i}
 $
is empty.
\end{itemize}
\end{Dn}
The main statements of this section are  Theorems \ref{partIsomThm} and \ref{complete goal} and
the most important technical result  is Lemma \ref{2.6}.   Among the technical instruments of the proof
of this Lemma  is the next Lemma \ref{a} which is useful in its own right.
\begin{La}
\label{a}
{\rm (\cite{Anton_Lebed}, Lemma 12.15).} Let $B$ be a $C^*$-subalgebra of the algebra
$L(H)$ of linear bounded operators in a Hilbert space $H$. If \/ $P_1, \  P_2\  \in\  B^\prime$
are two orthogonal projections such that the restrictions
$$
\left. B \right|_{H_{P_1}} \ \ and \ \  \left. B \right|_{H_{P_2}}
$$
(where $H_{P_1} = P_1 (H), \ \ H_{P_2} = P_2 (H)$) are both irreducible and
these restrictions are distinct representations then
$$
H_{P_1} \perp  H_{P_2} .
$$
\end{La}
\begin{La}
\label{2.6} Let $\V$ be a finely representable action such that
the corresponding complete interaction  $(\V,\H)$ be topologically
free. Let  $(C,\sigma,U)$ is a covariant representation of $\V$,
and  let $b$ be an  operator   of the form
\begin{equation}
\label{e--}
b= \sum_{x\in F} \sigma(a^{-x}_l)U_{x}^*\sigma(a^{-x}_r) + \sigma(a_0) + \sum_{x\in F} \sigma(a^{x}_l)U_{x}\sigma(a^{x}_r )
\end{equation}
 where $F $ is a finite subset of $\Gamma^+\setminus\{0\}$.
Then for every $\varepsilon >0$ there exists an irreducible representation
$\pi:\sigma (\A)\to L(H_{\pi })$ such that for any  irreducible
 representation $\nu: (\sigma\times U)(\A \rtimes_{(\V,\H)} \Gamma)\to L(H_\nu)$
 which is an extension of $\pi$  ($H_{\pi}\subset  H_\nu$) we have

(i)\ \ $\Vert \pi [\sigma (a_0)]\Vert \ge \Vert a_0 \Vert  - \varepsilon$,

(ii)\ \ $  P_{\pi } \,   \pi [\sigma (a_0)]\, P_{\pi } = P_{\pi } \, \nu (b)\,  P_{\pi }  $

\noindent where $P_{\pi }\in L(H_\nu)$ is the orthogonal projection onto $H_{\pi }$.
\end{La}
{\bf Proof.}\ \  As $\sigma$ is faithful we may and we shall
identify throughout the proof $\sigma (\A)$ and $\A$. For any
$a\in \A$ and $I\in \Prim\A$ we denote by $\breve{a} (I) $ the
number
\begin{equation}
\label{e2.1}
\breve{a} (I)  = \inf_{j\in I}\Vert  a+j \Vert
\end{equation}
For every $a\in \A$ the function $\breve{a} (\cdot) $ is lower semicontinuous on $\Prim \A$
and attains its upper bound equal to $\Vert a \Vert$ (see \cite{Dixmier}, 3.3.2. and 3.3.6.).
\\
Let $I_0 \in \Prim \A$ be a point at which  $\breve{a}_0 (I_0) =  \Vert a_0 \Vert$
and $\pi_0$ be an irreducible representation of $\A$ such that $I_0 = \ker \pi_0$ (thus $\Vert \pi_0 (a_0)\Vert = \Vert a_0 \Vert$).
 Since the function $\breve{a}_0 (\cdot)$ is lower semicontinuous it follows that for any $\varepsilon >0$
there exists an open set
$U \subset \Prim \A$ such that
\begin{equation}
\label{e2.1a}
\breve{a}_0 (I) >  \Vert a_0 \Vert - \varepsilon \ \ {\rm for\ \  every}\ \ I\in U.
\end{equation}
 As $F=\{x_1,...,x_k\}$ is finite and
 the interaction  $(\V,\H)$ is topologically free there exists a point
$I \in U $ such that all the points $t_{x_i}(I )$,  $i=1,...,k$ are   distinct (if they are
defined, i.e. if $I \in \Prim \A_{-x_i}$).
\\
Let $\pi$ be an
irreducible representation of $\A$ such that ${\rm ker}\, \pi = I$ and let $\nu$ be any extension of
$\pi$ up to an irreducible representation of   $(\pi\times U)(\A \rtimes_{(\V,\H)} \Gamma)$.
 For this representation $\nu$ we have
$$
H_{\pi}\subset H_\nu
$$
where $H_{\pi}$ is the representation space for $\pi$ and $H_\nu$
is that for $\nu$. Furthermore  for the orthogonal projection
$P_\pi:H_\nu\to H_\pi$ we have $\P_\pi\in \nu(\A)'$.
\\
By the choice of $\pi$ and (\ref{e2.1a}) we conclude that there exists a vector  $\xi \in H_{\pi}$
such that $\Vert \xi \Vert =1$ and
\begin{equation}
\label{e00}
\Vert {\pi} (a_0)\xi \Vert  > \Vert a_0\Vert -\varepsilon.
\end{equation}
Thus (i) is proved. \\
To prove (ii) it is enough to show  that for any monomials $\sigma(a^{x}_l)U_{x}\sigma(a^{x}_r)$,
$\sigma(a^{-x}_l)U_{x}^*\sigma(a^{-x}_r)$  which are elements of the sum \eqref{e2.1a} we have
\begin{equation}
\label{e2.2a}
P_{\pi}\,  \nu\Big(\sigma(a^{x}_l)U_{x}\sigma(a^{x}_r)\Big) \, P_{\pi} =0, \qquad P_{\pi}\,
 \nu\Big(\sigma(a^{-x}_l)U_{x}^*\sigma(a^{-x}_r)\Big) \, P_{\pi} =0.
\end{equation}
We will only prove the former relation as the proof for the latter one is completely analogous.
We fix an element $x$ in the set $F=\{x_1,...,x_k\}$ and consider the following possible positions of $I$.
\par
If $I \notin \Prim \A_{x}$ then  we have $ \nu (U_{x}U_{x}^*)P_{\pi}={\pi}(\V_{x}(1))=0$  and thus
$$
P_\pi \nu \big(\sigma(a^{x}_l)U_{x}\sigma(a^{x}_r)\big)P_\pi=P_\pi \nu \big(\sigma(a^{x}_l)U_{x}U^*_xU_x\sigma(a^{x}_r)\big)P_\pi
$$
$$
= \nu(\sigma(a^{x}_l))\nu (U_{x}U_{x}^*)P_\pi\nu \big(U_{x}\sigma(a^{x}_r)\big)P_\pi = 0.
$$
\par
If $I \notin \Prim \A_{-x}$ then observing that $ \nu(U_{x}^*U_{x})P_{\pi}={\pi}(\H_{x}(1))=0$
 we have
 $$
  P_\pi\nu \big(\sigma(a^{x}_l)U_{x}\sigma(a^{x}_r)\big)P_\pi= P_\pi
  \nu \big(\sigma(a^{x}_l)U_{x}\big)\nu (U_{x}^*U_{x})P_\pi\nu(\sigma(a^{x}_r))=0.
$$
\par
Finally let $I \in \Prim \A_{x} \cap \Prim A_{-x}$.\\
 In this case $\pi$
is an irreducible representation as for $\A_x$ so also for $\A_{-x}$ and
$t_{\pm x} (I)\in \Prim \A_x$ (according to the definition of $t_g$). Moreover  we have
\begin{equation}
\label{e2.3}
\nu (U_x^* U_x)\eta = \eta , \ \ \  \nu (U_x U_x^* )\eta = \eta  \ \ {\rm for \ \ any}
\ \ \eta \in H_{\pi}.
\end{equation}
In other words $H_{\pi}$ belongs as to the initial and final subspaces of $\nu (U_x)$
so also to the initial and final subspaces of $\nu (U_g^*)$.
\\
We will use Lemma \ref{a} where $P_1=P_{\pi}$ and $P_2 = \nu (U_x)P_{\pi}\nu (U_x^*)$.
By the definition of $\nu$ we have
$P_1\in \nu (\A)^\prime$  and (\ref{e2.3}) means that $
P_1 = P_1\, \nu (U_x^* U_x) = P_1\,  \nu (U_x U_x^* ) .
$ Moreover
$$
\nu (U_x): P_1(H_{\nu}) \to P_2(H_{\nu})
$$
 is an isomorphism. Observe also that
\begin{equation}
\label{e2.6}
P_2 \in \nu (\A_x)^\prime .
\end{equation}
Indeed, for any $a\in \A_x$ we have
$$
P_2 \nu (a)  =\nu (U_x)P_1 \nu (U_x^*)\nu (a) =  \nu (U_x)P_1 \nu (U_x^*) \nu (U_x U_x^*)\nu (a) =
$$
$$
\nu (U_x)P_1 \Big(\nu (U_x^*)\nu (a) \nu (U_x)\Big) \nu (U_x^*) =  \nu (U_x)P_1 \nu \big(\H_{x}(a)\big)\nu (U_x^*)=
$$
$$
 \nu (U_x)\nu (\H_{x}(a))\nu (U_x^*U_x)P_1 \nu (U_x^*)=
\nu \big(\V_x (\H_{x}(a))\big)\nu (U_x)P_1\nu (U_x^*)=\nu (a)P_2.
$$
Thus (\ref{e2.6}) is true.
In addition the irreducibility of
$\left. {\nu (\A_x)} \right|_{H_{P_1}}$ implies the irreducibility of
$\left. {\nu (\A_x)} \right|_{H_{P_2}}$ (here $H_{P_1} = P_1 (H_\nu )=H_{\pi}$ and
$H_{P_2} = P_2 (H_\nu )$).
\\
Now observe that for $a\in \A_x$ we have
$$
P_1  \nu (a)=0 \Longleftrightarrow \pi(a)P_1 = 0  \Longleftrightarrow \breve{a}(I )=0
$$
and
$$
P_2  \nu (a) =0 \Longleftrightarrow  \nu (U_x) P_1 \nu (U_x^*)\nu (a\V_x(1)) = 0
\Longleftrightarrow
$$
$$
\nu (U_x) P_1  \nu \left( \H_{x} (a)\right)\nu(U_x^*) = 0
\Longleftrightarrow
\nu (U_x)P_1  \nu \left( \H_{x} (a)\right) P_1 \nu(U_x^*) = 0
\Longleftrightarrow
$$
$$
P_1  \nu \left( \H_{x} (a)\right) =0 \Longleftrightarrow  \pi ( \H_{x} (a)) =0 \Longleftrightarrow
\breve{\H_{x}}(a) (I)=0\Longleftrightarrow \breve{a} (t_g (I)) =0.
$$
Thus, since the points $I$ and $t_x (I)$ are distinct we conclude that the representations
$\left. \nu (\A_x)\right|_{H_{P_1}}$ and $\left. \nu (\A_x) \right|_{H_{P_2}}$ are distinct. Applying Lemma \ref{a}
we find that
$$
P_1\cdot P_2 =0
$$
from which we have
$$
P_{\pi}\nu (U_x)P_{\pi}=P_{\pi}\nu (U_x)\nu (U_x^*U_x)P_{\pi}=P_{\pi}\nu (U_x)P_{\pi}\nu (U_x^*U_x)=P_1P_2 \nu (U_x)=0.
$$
Thus
$$
P_{\pi}\,  \nu\Big(\sigma(a^{x}_l)U_{x}\sigma(a^{x}_r)\Big) \, P_{\pi}=\nu(\sigma(a^{x}_l))\big(P_{\pi}
\nu(U_{x})  P_{\pi}\big)\nu(\sigma(a^{x}_r))=0
$$
which finishes the proof of  \eqref{e2.2a} and therefore the proof of the theorem as well.
\qed
\\
As a consequence, in the presence of topological freedom, we get that all covariant
representations satisfy a 'weaker version of $(^*)$ property'.
\begin{Tm}\label{partIsomThm}
Let $\V$ be a finely representable action such that the
corresponding complete interaction  $(\V,\H)$ be topologically
free. Then for every element $a\in C_0$     of the form
\begin{equation}\label{wertyuiop}
a= \sum_{x\in F} a^{-x}_l\hat{U}_{x}^*a^{-x}_r + a_0 + \sum_{x\in F} a^{x}_l\hat{U}_{x} a^{x}_r
\end{equation}
 where $F $ is a finite subset of $\Gamma^+\setminus\{0\}$ and for every $(C,\sigma,U)$
  covariant representation of $\V$, the  operator $(\sigma\times U)(a)$ determines uniquely
   the coefficient $a_0$. Namely for every $a$ of the form \eqref{wertyuiop} the  following inequality holds
$$
\|E_0(a)\| \leq \|(\sigma\times U)(a)\|.
$$

\end{Tm}
Since  in the case of covariance algebras   for $C^*$-dynamical systems, see Remark \ref{borek jest wolny},
 each finite sum of quasi-monomials may be presented in the form \eqref{wertyuiop} we immediately obtain the
 following statement, cf. also Proposition \ref{complete transfer action}.
\begin{Tm}\emph{\textbf{(Isomorphism Theorem for $C^*$-dynamical systems)}}\label{complete goal}\\
Let  $(\V,\H)$ be a topologically free complete  interaction  such that the hypotheses of Theorem
\ref{C-dynamical conditions}  hold.
 Then $\V$ and $\H$ are finely representable and for any covariant representation $(C,\sigma ,U)$ the formulae
$$
(\sigma\times U)(a)=\sigma(a), \qquad (\sigma\times U)(\hat{U}_x)=U_x,\quad \quad a\in\A, \,\,x\in \Gamma^{+}.
$$
determines the isomorphism  $(\sigma\times U)$  from  $\A \rtimes_{(\V,\H)} \Gamma$ onto
$(\sigma\times U)(\A \rtimes_{(\V,\H)} \Gamma)$.
\end{Tm}
In view of the foregoing statement Theorem \ref{partIsomThm} may be regarded as a part of
Isomorphism Theorem for
interactions. However, the following problem still remains open.
\begin{quote}
\textbf{Problem.} Can the Isomorphism Theorem be extended  to the general finely representable actions case?
\end{quote}

\end{document}